\documentclass{amsart} 

\usepackage{amsmath,amsthm}     
\usepackage{amssymb}            
\usepackage{euscript}           
\usepackage{enumerate,calc}
\usepackage[matrix,arrow,curve,frame]{xy}    

\xymatrixcolsep{1.9pc}                          
\xymatrixrowsep{1.9pc}
\newdir{ >}{{}*!/-5pt/\dir{>}}                  

\newtheorem{thm}[subsection]{Theorem}
\newtheorem{defn}[subsection]{Definition}
\newtheorem{prop}[subsection]{Proposition}

\newtheorem{cor}[subsection]{Corollary}
\newtheorem{lemma}[subsection]{Lemma}

\theoremstyle{definition}  

\newtheorem{remark}[subsection]{Remark}

  {\end{list}}
  
\newcommand{\dfn}{\textbf} 

\newcommand{\mdfn}[1]{\dfn{\mathversion{bold}#1}} 

\newcommand{\tens}              {\otimes}               
\newcommand{\iso}               {\cong}  

\newcommand{\cat}{\EuScript}    
\newcommand{\cC}{{\cat C}}

\newcommand{\field}[1]  {\mathbb #1} 
\newcommand{\A}         {\field A}
\newcommand{\R}         {\field R}

\newcommand{\M}         {\field M}

\newcommand{\Z}         {\field Z}
\newcommand{\C}         {\field C}
\renewcommand{\P}         {\field P}

\DeclareMathOperator{\spec}{Spec}
\DeclareMathOperator{\Spec}{Spec}

\DeclareMathOperator{\CH}{CH}
\DeclareMathOperator{\coker}{coker}
\DeclareMathOperator{\chara}{char}

\newcommand{\ra}{\rightarrow}                   
\newcommand{\lra}{\longrightarrow}              
\newcommand{\la}{\leftarrow}                    
\newcommand{\lla}{\longleftarrow}               
\newcommand{\llra}[1]{\stackrel{#1}{\lra}}      
\newcommand{\llla}[1]{\stackrel{#1}{\lla}}      

\newcommand{\inc}{\hookrightarrow}              

\newcommand{\blank}{-}                          
\newcommand{\Id}{Id}                            

\newcommand{\pt}{pt}

\newcommand{\rup}[1]{\lceil \frac{#1}{2} \rceil}

\newcommand{\rea}[1]{|{#1}|}             
\newcommand{\map}{\rightarrow}

\newcommand{\ceck}[1]{\Cech(#1)}         
\newcommand{\oceck}[1]{\Cech^{o}(#1)}    
\newcommand{\oreal}[1]{\rea{\oceck{U}}}  
\newcommand{\creal}[1]{\rea{\ceck{U}}}   

\newcommand{\Cech}{\check{C}}

\newcommand{\tH}{\tilde{H}}

\newcommand{\RP}{\R{\text{\sl P}}}

\newcommand{\cone}{/\!\!\!\!\Sigma}

\newcommand{\DQ}{DQ}

\numberwithin{equation}{subsection}

\newenvironment{myequation}
  {\addtocounter{subsection}{1}\begin{eqnarray}}
  {\end{eqnarray}$\!\!$}

\newcounter{property}

\newenvironment{property}
  {\addtocounter{property}{1}
   \medskip
   \noindent {\bf Property \Alph{property}.}~}
  {\medskip}

\begin{document}

\title{The Hopf condition for bilinear forms
over arbitrary fields}

\author{Daniel Dugger}
\author{Daniel C. Isaksen}

\address{Department of Mathematics\\ University of Oregon\\ Eugene, OR
97403}

\address{Department of Mathematics\\ University of Notre Dame\\ 
Notre Dame, IN 46556}

\email{ddugger@math.uoregon.edu}

\email{isaksen@math.wayne.edu}

\begin{abstract}
We settle an old question about the existence of certain
`sums-of-squares' formulas over a field $F$, related to the
composition problem for quadratic forms.  A classical theorem says
that if such a formula exists over a field of characteristic $0$, then
certain binomial coefficients must vanish.  We prove that this result
also holds over fields of characteristic $p>2$.
\end{abstract}

\maketitle


\section{Introduction}
\label{se:intro}

Fix a field $F$.  A classical problem asks for what values of $r$,
$s$, and $n$ do there exist identities of the form
\begin{myequation}
\label{eq:main}
 \Bigr ( \sum_{i=1}^{r}x_i^2 \Bigl ) \cdot
  \Bigr ( \sum_{i=1}^{s}y_i^2 \Bigl ) =
  \sum_{i=1}^{n}z_i^2 
\end{myequation}
where the $z_i$'s are bilinear expressions in the $x$'s and $y$'s.
Equation (\ref{eq:main}) is to be interpreted as a formula in the
polynomial ring $F[x_1,\ldots,x_r,y_1,\ldots,y_s]$; we call it a
\mdfn{sum-of-squares formula of type $[r,s,n]$}.

The question of when such formulas exist has been extensively studied:
\cite{La} and \cite{S1} are excellent survey articles, and \cite{S2}
is a detailed sourcebook.  In this paper we prove the following
result, solving Problem C of \cite{La}:

\begin{thm}
\label{th:main}
If $F$ is a field of characteristic not equal to $2$, and a
sum-of-squares formula of type $[r,s,n]$ exists over $F$, then
$\binom{n}{i}$ must be even for $n-r<i<s$.
\end{thm}

We now give a little history.  It is common to let $r*_F s$ denote the
smallest $n$ for which a sum-of-squares formula of type $[r,s,n]$
exists.   Many papers have studied lower
bounds on $r *_F s$, but for a long time such results  were known only
for fields of characteristic $0$: one reduces to a geometric
problem over $\R$, and then topological methods are used to obtain the
bounds (see \cite{La} for a summary).  In this paper we begin the
process of extending such results to characteristic $p$, 
replacing the topological methods by those of motivic homotopy theory.

The most classical result along these lines is Theorem~\ref{th:main}
for the particular case $F=\R$, which leads to lower bounds
for $r*_\R s$.  It seems to have been proven in three places, namely
\cite{B}, \cite{Ho}, and \cite{St}; but in modern times the given
condition on binomial coefficients is usually called the `Hopf
condition'.  The paper \cite{S1} gives some history, and explains how
K.~Y. Lam and T.~Y. Lam deduced the condition for arbitrary fields of
characteristic $0$.  Problem C of \cite[p. 188]{La} explicitly asked whether
the same condition holds over fields of characteristic $p>2$.  Work on
this question had previously been done by Adem \cite{A1,A2} and
Yuzvinsky \cite{Y} for special values of $r$, $s$, and $n$.  In
\cite{SS} a weaker version of the condition was proved for arbitrary
fields and arbitrary values of $r$, $s$, and $n$.

Stiefel's proof of the condition for $F=\R$ used Stiefel-Whitney
classes; Behrend's (which worked over any formally real field) used
some basic intersection theory; and Hopf deduced it using singular
cohomology.  Our proof of the general theorem uses a variation of
Hopf's method and {\it motivic\/} cohomology.  It can be regarded as
purely algebraic---at least, as `algebraic' as things like group
cohomology and algebraic $K$-theory.  These days it is perhaps not so
clear that there exists a point where topology ends and algebra
begins.

\medskip

We now explain Hopf's proof, and our generalization, in more detail.
Given a sum-of-squares formula of type $[r,s,n]$, one has in
particular a bilinear map $\phi\colon F^r \times F^s \ra F^n$
given by $(z_1, \ldots, z_n)$.
If we let $q$ be the quadratic form on $F^k$ given by
$q(w_1,\ldots,w_k)=w_1^2+\cdots+w_k^2$, then we have
$q(\phi(x,y))=q(x) q(y)$.  When
$F=\R$ one has that $q(w)=0$ only when
$w=0$, and so $\phi$ restricts to a map $(\R^r-0)\times
(\R^s-0) \ra (\R^n-0)$.  The bilinearity of $\phi$ tells us in
particular that we can quotient by scalar-multiplication to get
$\RP^{r-1} \times \RP^{s-1} \ra \RP^{n-1}$.

On mod 2 cohomology this gives $\Z/2[x]/x^n \ra \Z/2[a]/a^r \tens
\Z/2[b]/b^s$, and the bilinearity of $\phi$ shows that $x\mapsto a+b$.
Since $x^n=0$ and we have a ring map, it follows that $(a+b)^n=0$ in
the target ring.  The Hopf condition falls out immediately.

This proof used, in a seemingly crucial way, the fact that over $\R$ a
sum of squares is $0$ only when all the numbers were zero to begin
with.  This of course does not work over fields of characteristic $p$
(or over $\C$, for that matter).  Our bilinear form gives us a map of
schemes $\phi\colon\A^r \times \A^s \ra \A^n$, but we cannot say that
it restricts to $(\A^r-0)\times (\A^s-0) \ra (\A^n-0)$ as we did
above.

To remedy the situation, let $Q_k$ denote the projective quadric in
$\P^{k+1}$ defined by the equation $w_1^2+\cdots+w_{k+2}^2=0$.  The
bilinear map $\phi$ induces
\[ (\P^{r-1}-Q_{r-2}) \times (\P^{s-1}-Q_{s-2}) \ra (\P^{n-1}-Q_{n-2}).
\]
In effect, we have removed all possible numbers whose sum-of-squares
would give us zero.  Let $DQ_k$ denote the \dfn{deleted quadric}
$\P^k-Q_{k-1}$ (our convention is that the subscript on a scheme
always denotes its dimension).  We will compute the mod $2$ motivic
cohomology of $DQ_k$ (Theorem~\ref{th:DQcohom}), find that it is close
to being a truncated polynomial algebra, and repeat Hopf's argument in
this new context.  As an amusing exercise (cf. \cite[6.3]{Ln}) one can
show that over the field $\C$ the space $DQ_k$---with the complex
topology---has the same homotopy type as $\RP^{k}$; so our argument is
in some sense `the same' as Hopf's in this case.

The idea of using deleted quadrics to deduce the Hopf condition first
appeared in \cite{SS}.  In that paper the Chow groups of the deleted
quadrics were computed, but these are only enough to deduce a weaker
version of the Hopf condition (one that is approximately half as
powerful).  This is explained further in Remark~\ref{re:SS}.  On the
other hand, we should point out that the full power of motivic
cohomology is not completely necessary in this paper: one can also
derive the Hopf condition using \'etale cohomology, by the same
arguments (see Remark \ref{rem:etale}).  Since in this case computing
\'etale cohomology involves exactly the same steps as computing
motivic cohomology, we have gone ahead and computed the stronger
invariant (to be used in future work).  

\subsection{Organization}
Section \ref{se:basic} 
shows how to deduce the Hopf condition from a few easily
stated facts about motivic cohomology.  Section \ref{se:review}
outlines in more
detail the basic properties of motivic cohomology needed in the rest
of the paper.  This list is somewhat extensive, but our hope is that
it will be accessible to readers not yet acquainted with the motivic
theory---most of the properties are analogs of familiar things about
singular cohomology.  Finally, Section \ref{se:computation}
carries out the necessary calculations.
We also include an appendix on  the Chow groups of
quadrics, as several facts about these play a large role in the
paper.


\section{The basic argument}
\label{se:basic}

Because of the nature of the computations that we will make, we use
slightly different definitions for the varieties $Q_n$ and $DQ_n$ than
those in Section \ref{se:intro}.  These definitions will remain in
effect for the entire paper.  Unfortunately, the usefulness of these
choices will not become clear until Section \ref{se:computation}.

From now on the field $F$ is always assumed not to have characteristic
$2$.

\begin{defn}
When $n = 2k$, let $Q_n$ be the projective quadric in $\P^{n+1}$
defined by the equation $a_1 b_1 + a_2 b_2 + \cdots + a_{k+1} b_{k+1}=0$.
When $n = 2k+1$, let $Q_n$ be the projective quadric in $\P^{n+1}$
defined by the equation $a_1 b_1 + a_2 b_2 + \cdots a_{k+1} b_{k+1} + c^2=0$.
In either case, let $DQ_{n+1}$ be $\P^{n+1} - Q_n$.
\end{defn}

Note that $Q_0$ is isomorphic to $\spec F \amalg \spec F$, and $Q_1
\iso \P^1$.  One possible isomorphism $\P^1 \map Q_1$ sends $[x,y]$ to
$[-x^2, y^2, xy]$.

Occasionally we will need to equip $DQ_{n+1}$ with a basepoint, in
which case we'll always choose $[1,1,0,0,\ldots,0]$ (although the
choice turns out not to matter).

\begin{lemma}
\label{lem:coord-change}
Suppose that the ground field $F$ has a square root of $-1$ (call it $i$).
Then $Q_n$ is isomorphic to the projective quadric in $\P^{n+1}$
defined by the equation
$w_1^2 + \cdots + w_{n+2}^2 = 0$.  
\end{lemma}

\begin{proof}
When $n = 2k$, use the change of coordinates 
$a_j=w_{2j-1}+iw_{2j}$, $b_j=w_{2j-1}-iw_{2j}$.  
When $n = 2k+1$, use the
same formulas as above for $1 \leq j \leq k+1$ and also let $c=w_{n+2}$.
\end{proof}

We regard $\P^{2k} \inc \P^{2k+1}$ as the subscheme defined by
$a_{k+1} = b_{k+1}$, and we regard $\P^{2k-1} \inc \P^{2k}$ as the
subscheme defined by $c = 0$.  These choices have the advantage that
they give us inclusions $Q_{n-2} \inc Q_{n-1}$ and $DQ_{n-1} \inc
DQ_{n}$.

The following theorem states the computation of the motivic cohomology
ring $H^{*,*}(DQ_n;\Z/2)$. In order to understand the statement, the
reader needs to know just a few basic facts about motivic cohomology;
a more complete account of these facts appears in Section
\ref{se:review}.  First, $H^{*,*}(\blank;\Z/2)$ is a contravariant
functor defined on smooth $F$-schemes, taking its values in bi-graded
commutative rings of characteristic $2$.  If we set $\M_2 =
H^{*,*}(\spec F;\Z/2)$, the map induced by $X\ra \spec F$ makes
$H^{*,*}(X;\Z/2)$ into an $\M_2$-algebra.  It is known that
$\M_2^{0,0}\iso \Z/2$, $\M_2^{0,1}\iso \Z/2$, and the generator $\tau
\in \M_2^{0,1}$ is not nilpotent.

\begin{thm}
\label{th:DQcohom}
Assume that every element of $F$ is a square and that $\chara(F) \neq 2$.
\begin{enumerate}[(a)]
\item If $n=2k+1$ then $H^{*,*}(DQ_n;\Z/2) \iso \M_2[a,b]/(a^2=\tau
b, b^{k+1})$, where $a$ has degree $(1,1)$ and $b$ has degree $(2,1)$.
\item If $n=2k$ then $H^{*,*}(DQ_n;\Z/2) \iso \M_2[a,b]/(a^2=\tau b,
b^{k+1},ab^k)$ where $a$ and $b$ are as in part (a).
\item The
map $H^{*,*}(DQ_{n+1};\Z/2) \ra H^{*,*}(DQ_n;\Z/2)$ sends $a$ to $a$,
and $b$ to $b$.
\end{enumerate}
\end{thm}

In fact, $b$ is the unique nonzero class in $H^{2,1}$, and
$a$ is the unique nonzero class in $H^{1,1}$ which becomes zero when
restricted to the basepoint $\spec F \ra DQ_n$.  These facts are
needed below in the proof of Proposition \ref{pr:p^*}.  See
the comments before Proposition \ref{pr:DQmap} for more details.

Note that if $\tau$ were equal to $1$ then the above rings would be
truncated polynomial algebras (in analogy with the singular cohomology
of $\RP^n$).  

A more general version of this theorem, without any assumptions on
$F$, appears as Theorem~\ref{th:hyperbolic}.  The proof is slightly
involved, and so deferred until Section \ref{se:computation}.  
However, let us at least
record how the above statements follow from the more general version:

\begin{proof}
If every element of $F$ is a square, then $\M^{1,1}_2 = 0$ (see
Section \ref{se:finiteprop}).  Therefore, in Theorem
\ref{th:hyperbolic} both $\rho$ and $\epsilon$ are zero.  This gives
us the formulas in part (a) and (b).  Part (c) is Proposition
\ref{pr:DQmap}.
\end{proof}

For us, the most important consequence of the theorem is 
the following:

\begin{cor}
In $H^{*,*}(DQ_n;\Z/2)$ we have $a^{n+1}=0$ and $a^i \neq 0$ for
$i\leq n$.
\end{cor}

\begin{proof}
The claims are immediate from the calculation since all the
powers of $\tau$ are nonzero.
\end{proof}

\begin{proof}[Proof of Theorem~\ref{th:main}]
Suppose we have a sum-of-squares formula of type $[r,s,n]$ over $F$.
This remains true if we extend $F$, so we may as well assume 
that every element of $F$ is a square.
Therefore,
Theorem \ref{th:DQcohom} applies.

As explained in Section \ref{se:intro}, the sum-of-squares formula gives a map
$p\colon DQ^{r-1} \times DQ^{s-1} \ra DQ^{n-1}$ (this uses Lemma
\ref{lem:coord-change}) and we will consider the induced map on
motivic cohomology.  There is a K\"unneth formula for computing
motivic cohomology of products of certain `cellular' varieties (see
Proposition \ref{pr:kunneth}), and the deleted quadrics belong to this
class by Proposition \ref{pr:cell}.  In order to apply 
Proposition \ref{pr:kunneth}, we also have to observe that 
$H^{*,*}(DQ_{r-1};\Z/2)$ is free over $\M_2$, which is apparent from
Theorem \ref{th:DQcohom}.

Therefore
$p^*$ is a map
\[ H^{*,*}(DQ_{n-1};\Z/2) \ra H^{*,*}(DQ_{r-1};\Z/2) \tens_{\M_2}
H^{*,*}(DQ_{s-1};\Z/2).
\]
We will use the letters $a$ and $b$ to denote the generators of
$H^{*,*}(DQ_{n-1};\Z/2)$, $a_1$ and $b_1$ for the generators of
$H^{*,*}(DQ_{r-1};\Z/2)$, and $a_2$ and $b_2$ for the generators of
$H^{*,*}(DQ_{s-1};\Z/2)$.

We show in the following proposition that $p^*(a)=a_1+a_2$.  Since the
above corollary says that $a^{n}=0$, it will follow that
$(a_1+a_2)^{n}=0$.  Using the corollary again, this can only happen if
$\binom{n}{i}$ is even for $n-r<i<s$.
\end{proof}

\begin{prop}
\label{pr:p^*}
Suppose that $F$ is a field of characteristic not 2 
in which every element is a square.  If
$p^*$, $a$, $a_1$, and $a_2$ are as in the above proof, then
$p^*(a)=a_1+a_2$.
\end{prop}

Before we can give the proof, we need to state a few more properties
of motivic cohomology.  Once again, more details are given in Section
\ref{se:review}.  First, $\M_2^{p,q}$ is nonzero only in the range
$q\geq 0$.  Second, when every element of $F$ is a square one has
$\M^{1,1}_2 = 0$.  Finally, motivic cohomology is $\A^1$-homotopy
invariant in the following sense.  Let $i_0$ and $i_1$ denote the
inclusions $\{0\} \inc \A^1$ and $\{1\}\inc \A^1$, respectively.  If
$H\colon X\times \A^1 \ra Y$ is a map of smooth schemes, then the
composites $H(\Id \times i_0)$ and $H(\Id \times i_1)$ induce the same
map $H^{*,*}(Y;\Z/2) \ra H^{*,*}(X;\Z/2)$.  Such a map $H$ is called
an `$\A^1$-{\it homotopy\/}' from $H(\Id \times i_0)$ to $H(\Id \times
i_1)$.

\begin{proof}
Because $p^*(a)$ has degree $(1,1)$, it must be of the form
$\epsilon_1 a_1 + \epsilon_2 a_2 + m\cdot 1$, where $m$ belongs to
$\M^{1,1}_2$ and $\epsilon_1$ and $\epsilon_2$ belong to $\M_2^{0,0}\iso \Z/2$.
Since $\M^{1,1}_2=0$ under our assumptions on $F$, we can ignore $m$.
To show that $\epsilon_1 = 1$,
in light of Theorem \ref{th:DQcohom}(c) it would suffice to verify that the map
\[ DQ_1 \times \{*\} \ra DQ_{r-1}\times DQ_{s-1} \ra DQ_{n-1} 
\]
is $\A^1$-homotopic to the standard inclusion $DQ_1 \inc DQ_{n-1}$.
(A similar argument will show that $\epsilon_2=1$.)  Actually we will
not quite do this, but instead verify that the composition
\[ j\colon DQ_1 \times \{*\} \ra DQ_{r-1}\times DQ_{s-1} \ra DQ_{n-1}
\inc DQ_{n+1} 
\]
is $\A^1$-homotopic to the standard inclusion.
By Theorem \ref{th:DQcohom}(c) again, this is enough.

For the rest of this section we will use the coordinates
$w_1,\ldots,w_{n+2}$ on $\P^{n+1}$ given in
Lemma~\ref{lem:coord-change}.  Recall that $\phi$ is our bilinear map
$F^r \times F^s \ra F^n$.  Let ${e_1, \ldots, e_k}$ 
be the standard basis for $F^k$,
and let $\phi(e_1,e_1)=(u_1,\ldots,u_n)$ and
$\phi(e_2,e_1)=(v_1,\ldots,v_n)$.  Then the map $j\colon DQ_1 \ra DQ_{n+1}$
has the form
\[ [a,b] \mapsto [u_1a+v_1b,u_2a+v_2b,\ldots,u_na+v_nb,0,0],
\]
and the sum-of-squares formula satisfied by $\phi$ tells us
that 
\[ u_1^2+\cdots+u_n^2=1, \quad v_1^2+\cdots+v_n^2=1, \quad \text{and}
\quad
u_1v_1+\cdots+u_nv_n=0.
\]
Note that the standard inclusion $DQ_1 \inc DQ_{n+1}$ has the same
description but where $(u_1,\ldots,u_n)=(1,0,\ldots,0)$ and
$(v_1,\ldots,v_n)=(0,1,0,\ldots,0)$.  The following lemma gives the
desired $\A^1$-homotopy, since both the map $j$ and the standard
inclusion are homotopic to the map $[a,b] \mapsto [0,0,\ldots,0,a,b]$.
\end{proof}

For the following statement, recall that we are still using the coordinates on
$\P^{n+1}$ given by Lemma~\ref{lem:coord-change}.
 
\begin{lemma}
Suppose that $F$ contains a square root of $-1$.  Let $u$ and $v$ be
vectors in $F^n$ such that $\Sigma_j u_j^2 = 1 = \Sigma_j
v_j^2$ and $\Sigma_j u_j v_j = 0$.  Then the map $f:DQ_1 \map
DQ_{n+1}$ given by
\[ [a,b] \mapsto [u_1a+v_1b,u_2a+v_2b,\ldots,u_na+v_nb,0,0]
\]
is $\A^1$-homotopic to the map 
$[a,b] \mapsto [0,0,\ldots,0,a,b]$.
\end{lemma}

\begin{proof}
Let $i$ be a square root of $-1$.
First define a homotopy $DQ_1 \times \A^1 \map DQ_{n+1}$ by
the formula
\[
([a,b],t) \mapsto [u_1a+v_1b,u_2a+v_2b,\ldots,u_na+v_nb,ta-tib,tia+tb].
\]
This shows that $f$ is homotopic to $g$,
where $g$ is the map
\[
[a,b] \mapsto [u_1a+v_1b,u_2a+v_2b,\ldots,u_na+v_nb,a-ib,ia+b].
\]
Now define
another homotopy $DQ_1 \times \A^1 \map \P^{n+1}$ by
the formula
\[
([a,b],t) \mapsto [tu_1a+tv_1b,tu_2a+tv_2b,\ldots,tu_na+tv_nb,a-tib,tia+b].
\]
The assumptions on the $u$'s and $v$'s imply that the sum of the
squares in the image is exactly equal to $a^2+b^2$, which is nonzero
because $[a,b]$ lies in $DQ_1$.  So this is actually a homotopy $DQ_1 \times
\A^1 \ra DQ_{n+1}$, showing that $g$ is homotopic to the desired map.
\end{proof}

\begin{remark}
\label{re:SS}
In \cite{SS} a weaker version of the Hopf condition was obtained by
computing the Chow ring $\CH^*(DQ_n)$, which essentially corresponds
to the subring of $H^{*,*}(DQ_n;\Z/2)$ generated by $b$ (see Property
(A) in Section \ref{se:review}).  
This amounts to seeing about half of what motivic
cohomology sees.
\end{remark}

\begin{remark}
\label{rem:etale}
When $F$ has a square root of $-1$, a theorem of \cite{Lv} says
that the \'etale cohomology ring $H^{*}_{et}(DQ_n;\mu_2^{\tens *})$ is
isomorphic to $H^{*,*}(DQ_n;\Z/2)[\tau^{-1}] \iso H^{*,*}(DQ_n;\Z/2)
\otimes_{\M_2} \M_2[\tau^{-1}]$ (see Property (I) below).  Since
$H^{*,*}(DQ_n;\Z/2)$ is free over $\M_2$, this localization is
particularly simple: it is precisely a truncated polynomial algebra
$\M_2[\tau^{-1}][a]/a^{n+1}$.  So the Hopf condition could have been
proven using \'etale cohomology.
\end{remark}

\begin{remark}
When every element of $F$ is a square, it follows from the proof of
the Milnor conjecture \cite{V3} that $\M_2 \iso \Z/2[\tau]$.  We never
needed this, but it's useful to keep in mind.
\end{remark}


\section{Review of motivic cohomology}
\label{se:review}

The theory now called motivic cohomology was first developed in two
main places, namely \cite{Bl1} and \cite{VSF} (together with many
associated papers).  The paper \cite{V4} proved that the two
approaches give isomorphic theories.  Below we recall the basic
properties of motivic cohomology needed in the paper.  For
various reasons it is difficult to give simple references to
\cite{VSF} so most of our citations will be to \cite[Sec.~3]{SV}
and the lecture notes \cite{MVW}.

\subsection{Basic properties}
\label{pr:motivic}
For every field $F$, motivic cohomology is a contravariant functor
$H^{*,*}(\blank)$ from the category of smooth schemes of finite type
over $F$ to the category of bi-graded commutative rings.
Commutativity means that if $a\in H^{p,q}(X)$ and $b\in H^{s,t}(X)$
then $ab=(-1)^{ps}ba$.  For the basic construction we refer the reader
to \cite[Sec.~3]{SV} or \cite[Sec.~3]{MVW}.  The list of properties
below is far from complete, and in some cases we only give crude
versions of more interesting properties---but this is all we will need
in the present paper.

The scheme $\spec F$ will often be denoted by ``$\pt$'', and we denote
$H^{*,*}(\pt)$ by $\M$.  The ring $\M$ can be very complicated (and
is, in general, unknown).  The motivic cohomology of a scheme is
naturally a graded-commutative algebra over $\M$.

\begin{property}
\label{prop:Chow}
The graded subring $\oplus_n H^{2n,n}(X)$ is naturally
isomorphic to the Chow ring $\CH^*(X)$.
\cite[p.~268]{Bl1}, \cite[p.~4; Lect.~17]{MVW}.
\end{property}

In particular, $\M^{0,0} = \Z$.  In general, $H^{*,*}(X)$ is
isomorphic to the higher Chow groups of $X$ \cite[Cor.~1.2]{V4}.

\begin{property}
\label{prop:local}
For a closed inclusion $j\colon Z\inc X$ of smooth schemes of
codimension $c$, there is a long exact sequence of
the form
\[  \cdots  \ra H^{*-2c,*-c}(Z) \llra{j_!} H^{*,*}(X) \ra H^{*,*}(X-Z) 
              \ra H^{*-2c+1,*-c}(Z) \ra \cdots 
\]
The map $j_!$ is called the `Gysin map' or the `pushforward', and it
is a map of $\M$-modules.  The long exact sequence is called
the {\it Gysin\/}, {\it localization\/}, or {\it purity\/} sequence.
\cite[Sec.~3]{Bl1}, \cite{Bl2}.
\end{property}

\begin{property}
\label{prop:homotopy}
Let $i_0$ and $i_1$ denote the inclusions $\{0\} \inc \A^1$ and
$\{1\}\inc \A^1$, respectively.  
If $H\colon X\times \A^1 \ra Y$ is a map of smooth schemes,
then the composites $H(\Id \times i_0)$ and $H(\Id \times i_1)$ 
induce the same map $H^{*,*}(Y)
\ra H^{*,*}(X)$.  Such a map $H$ is called an $\A^1$-homotopy
from $H(\Id \times i_0)$ to $H(\Id \times i_1)$.
\cite[Sec.~2]{Bl1}, \cite[Prop.~4.2]{SV}.
\end{property}

\begin{property}
\label{prop:P}
$H^{*,*}(\P^n) = \M[t]/(t^{n+1})$, where $t$ has degree $(2,1)$. 
\cite[Prop.~4.4]{SV}.
\end{property}

\begin{property}
\label{prop:bundle}
If $E\ra B$ is an algebraic fiber bundle (i.e., a map which is locally
a product in the Zariksi topology) whose fiber is an
affine space $\A^n$, then $H^{*,*}(B) \ra H^{*,*}(E)$ is an
isomorphism.
\end{property}

Property (E) is easy to prove by inducting on the size of a
trivializing cover, and using the Mayer-Vietoris sequence
\cite[Prop.~4.1]{SV} together with Property (C).

\begin{property}
\label{prop:negative}
$\M^{p,q}=0$ if  $q<0$, if $p>q\geq 0$, or if $q=0$ and $p<0$.
\cite[p.~4; Th.~3.5]{MVW}
\end{property}

\begin{property}
\label{prop:1}
$\M^{1,1}=F^*$ and $\M^{0,1}=0$.
\cite[Th.~6.1]{Bl1}, \cite[p.~4,(2)]{MVW}.
\end{property}

\subsection{Finite coefficients}
\label{se:finiteprop}
For every $n \in \Z$ there is also a theory $H^{*,*}(\blank;\Z/n)$
which is related to $H^{*,*}(\blank)$ by a natural long exact
sequence of the form
\begin{myequation}
\label{eq:coeffs}
 \cdots \ra H^{*,*}(X) \llra{\times n} H^{*,*}(X) \ra H^{*,*}(X;\Z/n) \ra
H^{*+1,*}(X) \llra{\times n} \cdots
\end{myequation}
For the definition see \cite[Def.~3.4]{MVW}.  The theory satisfies the
analogs of Properties (B) through (F) above.

Let $\M_2$ denote $H^{*,*}(\pt;\Z/2)$.  Since $\M$ may contain
$2$-torsion, $\M_2$ is not necessarily the same as
$\M/(2)$---rather, there is a long exact sequence of the form
\[ \cdots \ra \M^{p,q} \llra{\times 2} \M^{p,q} \ra \M_2^{p,q} \ra 
\M^{p+1,q} \ra \cdots 
\]
This sequence, together with Property (F) and the fact that
$\M^{0,0}=\Z$, tells us that $\M^{0,0}_2 = \Z/2$.  Note that
$H^{*,*}(X;\Z/2)$ is naturally a commutative algebra over $\M_2$.

Since $\M^{1,1}=F^*$ and
$\M^{0,1}=\M^{2,1}=0$, we get
the exact sequence
\begin{myequation}
\label{eq:square}
 0 \ra \M_2^{0,1} \ra F^* \llra{\times 2} F^* \ra \M_2^{1,1} \ra 0
\end{myequation}
where the map $F^* \ra F^*$ sends $x$ to $x^2$.  The usual notation is
to let $\tau \in \M_2^{0,1}$ denote the class which maps to $-1$, and
to let $\rho\in \M_2^{1,1}$ denote the image of $-1$.  If $F$ has a
square root of $-1$ then $\rho=0$.  Moreover, if every element of $F$
is a square then $\M_2^{1,1} = 0$.

\subsection{The Bockstein}
\label{se:bockstein}
The Bockstein map $\beta\colon
H^{*,*}(\blank;\Z/2) \ra H^{*+1,*}(\blank;\Z/2)$ is defined in
the usual manner from the maps in the sequence (\ref{eq:coeffs}).
A direct consequence of the definition (as in topology) is that $\beta^2 = 0$.
Note that $\beta(\tau)=\rho$.

\enlargethispage{\baselineskip}

\begin{property}
\label{prop:Bockstein}
For all $a,b\in H^{*,*}(X;\Z/2)$, $\beta(ab)=\beta(a)b+a\beta(b)$.
\cite[Lem.~6.1]{Lv}.
\end{property}

\subsection{Relation with \'etale cohomology}
\label{se:etale}

There is a natural map of bi-graded rings $\eta\colon H^{*,*}(X;\Z/n) \ra
H^*_{et}(X;\mu_n^{\tens *})$ (cf. \cite[Th.~10.2]{MVW}, for example).
In the case $n=2$, the element $\tau$ maps to the class of $-1$ in
$H^0_{et}(\pt;\mu_2)\iso \{1,-1\}$, and multiplication by this class
is an isomorphism on \'etale cohomology.
Note in particular that this implies that the powers of $\tau$ are all
nonzero in $H^{0,*}(\pt;\Z/2)$.

\begin{property}
\label{prop:etale}
The induced map
$H^{*,*}(X;\Z/2)[\tau^{-1}] \ra H^{*}_{et}(X;\mu_2^{\tens *})$ is an
isomorphism for any smooth scheme $X$, provided that
$F$ has a square root of $-1$. 
\cite{Lv}.
\end{property}

The construction of the map $\eta$ from \cite{MVW} makes it clear that the
Bockstein on $H^{*,*}(\blank;\Z/2)$ (which can be regarded as induced
by the extension $0\ra \Z/2 \ra \Z/4 \ra \Z/2\ra 0$) is compatible
with the Bockstein on \'etale cohomology induced by $0 \ra \mu_2 \ra
\mu_4 \ra \mu_2 \ra 0$.  If the field contains a square root of $-1$
then we can identify $\mu_4$ with $\Z/4$, and of course $\mu_2$ with
$\Z/2$.  These observations will be used in the proof of
Theorem~\ref{th:hyperbolic}.  

\subsection{Reduced cohomology}
\label{se:reduced}
Given any basepoint of a scheme $X$ (i.e., a map $\pt \map X$),
the kernel of the induced map $H^{*,*}(X) \map H^{*,*}(\pt)$
is the reduced cohomology of $X$ and is denoted by $\tH^{*,*}(X)$.
A similar definition applies with $\Z/n$-coefficients.  The above
map has a splitting (induced by $X \map \pt$), and thus
$H^{*,*}(X) \cong \M \oplus \tH^{*,*}(X)$.  Similarly,
$H^{*,*}(X;\Z/2) \cong \M_2 \oplus \tH^{*,*}(X;\Z/2)$.

\subsection{A K\"unneth theorem}
\label{se:kunneth}
Let $\cC$ denote the smallest class of smooth schemes satisfying
the following properties:
\begin{enumerate}[(1)]
\item $\cC$ contains the affine spaces $\A^k$.
\item If $Z\inc X$ is a closed inclusion of smooth schemes and $\cC$
contains two of $X$, $Z$, and $X-Z$, then it also contains the third.
\item If $E\ra B$ is an algebraic fiber bundle whose fiber is an
affine space, then $E\in \cC$ if and only if $B \in \cC$.
\end{enumerate}

The following result is a modest generalization of \cite[Th.~4.5]{J}, 
and can be proven using the same techniques.  A complete
proof, for a more general class of schemes than $\cC$, is given in
\cite{DI2}.

\begin{prop}
\label{pr:kunneth}
Suppose $X$ and $Y$ are smooth schemes, with at least one of them
belonging to $\cC$.  
If either $H^{*,*}(X)$ or $H^{*,*}(Y)$ is
free as an $\M$-module, then one has a K\"unneth isomorphism of
bi-graded rings
\[ H^{*,*}(X) \tens_{\M} H^{*,*}(Y) \iso H^{*,*}(X\times Y).
\]
Similarly, if either $H^{*,*}(X;\Z/n)$ or $H^{*,*}(Y;\Z/n)$ is
free as an $H^{*,*}(\pt;\Z/n)$-module, then one has a K\"unneth isomorphism of
bi-graded rings
\[ H^{*,*}(X;\Z/n) \tens_{H^{*,*}(\pt;\Z/n)} H^{*,*}(Y;\Z/n) 
   \iso H^{*,*}(X\times Y;\Z/n).
\]
\end{prop}


\section{Computations}
\label{se:computation}

In this section $F$ is an arbitrary ground field not of characteristic
$2$.  We will study the quadrics $Q_n$ and $DQ_n$.  Note that in this
generality Lemma \ref{lem:coord-change} does not apply; therefore,
$Q_n$ and $DQ_n$ cannot necessarily be redefined in terms of sums of
squares.  We assume $\chara(F)\neq 2$ so that $Q_{n}$ is smooth for
all $n$, not just even $n$.

\begin{prop}
\label{pr:free}
If $n$ is odd, $H^{*,*}(Q_n)$ is a free module over $\M$ with
generators in degrees $(0,0),(2,1),(4,2),\ldots,(2n,n)$.  If $n$ is
even, $H^{*,*}(Q_n)$ is a free $\M$-module with generators in degrees
$(0,0), (2,1),\ldots,(2n,n)$ plus an extra generator in degree
$(n,\frac{n}{2})$.
\end{prop}

\begin{proof}
The proof is by induction.  The result for $Q_0$ is obvious, and
the result for $Q_1 \cong \P^1$ is 
Property (D).

Except for the base cases in the previous paragraph, the argument for
the odd and even cases is identical.  We give details only for the
even case, so let $n = 2k$.  Let $Z$ be the $(n-1)$-dimensional
subscheme defined by $a_1 = 0$, and let $U=Q_n-Z$.  Note that $Z$ is
singular (it's the projective cone on $Q_{n-2}$), and that $U\iso
\A^{n}$.  Let $Q'=Q_n-\{[0, 1, 0, 0, \ldots,0]\}$, and let $Z'=Z-\{[0,
1, 0, 0, \ldots, 0]\}$.  Then $Z' \inc Q'$ is a smooth pair, with
complement $\A^n$.  So the localization sequence for $Z'\inc Q'$ gives us an
isomorphism $\tH^{*,*}(Q')\iso H^{*-2,*-1}(Z')$.  The projection map
$Z' \ra Q_{n-2}$ which forgets the first two homogeneous coordinates
is a fiber bundle with fiber $\A^1$, hence $H^{*,*}(Z')\iso
H^{*,*}(Q_{n-2})$ by Property (E).

Taking the computations of the previous paragraph together, we conclude
that $H^{*,*}(Q') \iso H^{*-2,*-1}(Q_{n-2}) \oplus \M$.  By
induction, this is free over $\M$ with one generator in each degree
$(0,0), (2,1), \ldots, (2n-2,n-1)$ plus an extra generator in degree
$(n, \frac{n}{2})$.

Finally, we consider the localization sequence for $\{[0, 1, 0,\ldots,0]\}
\inc Q_n$.   This has the form
\[ 
\cdots \la
H^{*-2n+1,*-n}(\pt) \llla{\delta} H^{*,*}(Q') 
\la H^{*,*}(Q_n) \la H^{*-2n,*-n}(\pt) \la \cdots
\]
The generators for $H^{*,*}(Q')$ (as an
$\M$-module) must map to zero under $\delta$ for dimension
reasons.  It follows that 
\[ 
0 \la H^{*,*}(Q') \la H^{*,*}(Q_n) \la H^{*-2n,*-n}(\pt) \la 0
\]
a short exact sequence of
$\M$-modules, in which the outer terms are known to be free.  So the
middle term is a direct sum of the outer terms.  The right term
provides a generator of degree $(2n,n)$, and the left term provides
the rest of the generators.
\end{proof}

The above proof also shows the following:

\begin{prop}
\label{pr:cell}
The schemes $Q_n$ and $DQ_n$ belong to the class $\cC$ from Section
\ref{se:kunneth}.  
\end{prop}

\begin{proof}
If one knows by induction that
$Q_{n-2}$ belongs to $\cC$ then so do $Z'$, $Q'$, and $Q_{n}$, in that
order.  

The fact that projective spaces belong to the class $\cC$ is
trivial: one uses the standard algebraic cell decomposition
(cf. \cite[1.9.1]{F}).  Then since $Q_{n-1}$ and $\P^n$ are both in
$\cC$, so is $DQ_n$.    
\end{proof}

By the above result, in order to understand the ring structure on
$H^{*,*}(Q_n)$ it suffices just to understand the subring
$H^{2*,*}(Q_n)\iso \CH^*(Q_n)$, because the $\M$-algebra
generators lie in degrees $(2*,*)$.  The computation of this Chow ring
is well-known; the additive computation can be found in
\cite[13.3]{Sw}, for instance, and the ring structure is stated in
\cite{KM}.  For the reader's convenience, and because we need several
of the auxiliary facts, we give a complete account in Appendix A.
These ideas lead to the following result, whose proof is essentially
the content of Theorem \ref{th:CHoddquadric} and
Theorem \ref{th:CHevenquadric}.

\begin{prop}\mbox{}\par
\begin{enumerate}[(a)]
\item
If $n=2k+1$ then as a ring
$H^{*,*}(Q_n)=\M[x,y]/(x^{k+1}-2y,y^2)$  where $x$
has degree $(2,1)$ and $y$ has degree $(2k+2,k+1)$.  
\item If $n=2k$ and $k$ is odd,
then $H^{*,*}(Q_n)=\M[x,y]/(x^{k+1}-2xy,y^2)$
where $x$ has degree $(2,1)$
and $y$ has degree $(2k,k)$.
\item If $n=2k$ and $k$ is even, 
then $H^{*,*}(Q_n)=\M[x,y]/(x^{k+1}-2xy, x^{k+1}y, y^2-x^ky)$
where $x$ has degree $(2,1)$
and $y$ has degree $(2k,k)$.
\end{enumerate}
\end{prop}

We will now consider the motivic cohomology of the deleted quadrics
$DQ_n$.  The idea is to use the localization sequence
\addtocounter{subsection}{1}
\begin{equation}
\label{eq:purity}
{}
\xymatrixcolsep{1pc}
\xymatrix{
\cdots &  H^{*-1,*-1}(Q_{n-1})\ar[l] &
H^{*,*}(DQ_n) \ar[l]
& H^{*,*}(\P^n) \ar[l]_-{i^*}
&H^{*-2,*-1}(Q_{n-1}) \ar[l]_-{j_!}
& \cdots\ar[l].}
\end{equation}
By Proposition \ref{pr:free} the cohomology of $Q_{n-1}$ has
generators as an $\M$-module in degrees $(2*,*)$, so we can completely
determine
the $\M$-module map $j_!$ just by understanding the pushforward map
$\CH^{*-1}(Q_{n-1}) \map \CH^*(\P^{n})$ of Chow groups.  For the quadrics,
this is discussed in detail in the appendix: all maps are either the
identity or multiplication by $2$.
However, a problem now occurs.
Because the ground ring $\M$ might have $2$-torsion, the kernel and
cokernel of $j_!$ will not necessarily be free over $\M$---so we run
into complicated extension problems.  As a result, we haven't been
able to compute the integral motivic cohomology of $DQ_n$.  The
problem goes away if we work with $\Z/2$ coefficients.

\begin{prop}
\label{pr:hyperbolic-additive}
If $F$ is a field with $\chara(F) \neq 2$, then
$H^{*,*}(DQ_n;\Z/2)$ is a free $\M_2$-module with
one generator in degree 
$(i, \rup{i})$ for each $0 \leq i \leq n$, where $\rup{i}$ is the
smallest integer that is at least $\frac{i}{2}$.
\end{prop}

\begin{proof}
The argument from Proposition~\ref{pr:free} shows 
that $H^{*,*}(Q_{n-1};\Z/2)$ is free over $\M_2$ on the same set of
generators as before, and the map of subrings $H^{2*,*}(Q_{n-1}) \ra
H^{2*,*}(Q_{n-1};\Z/2)$ is just quotienting by the ideal $(2)$.

By Lemma \ref{lem:Gysin},
we know that the Gysin map
$j_!\colon H^{2i,i}(Q_{n-1}) \ra H^{2i+2,i+1}(\P^{n})$ is multiplication by 2
for $0 \leq i < \frac{n-1}{2}$, and is an isomorphism for 
$\frac{n-1}{2} < i \leq n-1$.  
If $n$ is odd, then it is the fold map
$\Z\oplus \Z \ra \Z$ for $i=\frac{n-1}{2}$.  

The goal is to use the $\Z/2$-analog of (\ref{eq:purity}), so
we first have to understand the Gysin map $j_!$ with $\Z/2$-coefficients.
Since $H^{2*,*}(Q_{n-1};\Z/2)$ and $H^{2*,*}(\P^n;\Z/2)$ are both
obtained from integral cohomology simply by quotienting by the ideal (2),
it follows that the Gysin map with $\Z/2$-coefficients is an
isomorphism, zero, or the fold map 
in all degrees $(2*,*)$.  Since the generators
(as $\M_2$-modules) live in these degrees, 
we find that the kernel
and cokernel of $j_!\colon H^{*,*}(Q_{n-1};\Z/2) \ra H^{*,*}(\P^n;\Z/2)$ are
both free over $\M_2$.  

If $n=2k$, then the generators for $\coker j_!$ are in degrees
$(0,0)$, $(2,1)$, $(4,2), \ldots,(2k,k)$, and the generators for $\ker
j_!$ are in degrees $(0,0), (2,1), \ldots, (2k-2,k-1)$.  If $n=2k+1$,
then the generators are the same, except that $\ker j_!$ has another
generator in degree $(2k,k)$.

From the $\Z/2$-analog of (\ref{eq:purity}), 
we have the short exact sequence
\[ 0 \la \ker j_! \la H^{*,*}(DQ_n;\Z/2) \la \coker j_! \la 0.
\]
It follows that the middle group is also free over $\M_2$.  
Be aware that the left map shifts degrees by $(-1,-1)$.  
\end{proof}

We know $\M_2^{0,0}\iso \Z/2$.  From Property (F)
it follows that
$\M_2^{p,q}=0$ if $q<0$, if $q=0$ and $p<0$, or if $p>q \geq 0$.  So
the above calculation shows that $H^{1,1}(DQ_n;\Z/2) \iso
\M_2^{1,1}\oplus \M_2^{0,0}$, where the first summand comes from the
motivic cohomology of $\Spec F$.  Hence, there is a unique
nonzero element $a\in \tH^{1,1}(DQ_n;\Z/2)$.  When $n> 1$ the
calculation gives $H^{2,1}(DQ_n;\Z/2)\iso \Z/2$, and we let $b$ denote
the unique nonzero element.  For $n=1$ we have $DQ_1 \iso \A^1-0$,
and it is known that $H^{2,1}(\A^1-0;\Z/2)=0$ (see, for instance,
\cite[Lem.~6.8]{V2}).  
In this case we define $b=0$ by convention.

\begin{prop}
\label{pr:DQmap}
The map $H^{*,*}(DQ_{n+1};\Z/2) \map H^{*,*}(DQ_n;\Z/2)$ induced
by the inclusion takes $a$ to $a$ and $b$ to $b$.
\end{prop}

\begin{proof}
In light of the definitions of $a$ and $b$ in the previous
paragraph,
we just need to show that 
$H^{i,1}(DQ_{n+1};\Z/2) \map H^{i,1}(DQ_n;\Z/2)$ is surjective
for $i = 1$ or $i=2$.
Consider the diagram 
\[
\xymatrixcolsep{0.7pc}
\xymatrix{
H^{i+1,1}(\P^{n+1};\Z/2) \ar[d] & H^{i-1,0}(Q_n;\Z/2) \ar[d]\ar[l] &
    H^{i,1}(DQ_{n+1};\Z/2) \ar[d]\ar[l] & H^{i,1}(\P^{n+1};\Z/2) \ar[d]\ar[l]\\
H^{i+1,1}(\P^{n};\Z/2) & H^{i-1,0}(Q_{n-1};\Z/2) \ar[l] &
    H^{i,1}(DQ_{n};\Z/2) \ar[l] & H^{i,1}(\P^{n};\Z/2) \ar[l]   }
\]
in which the rows are localization sequences.  
The map between the groups on the left is the
identity.  A similar remark applies to the groups on the right.
Finally, the cohomology groups of the quadrics are both isomorphic to
$\M_2^{i-1,0}$, and the map between them is also the identity.  It
follows from a diagram chase that the desired map is surjective.
\end{proof}

\begin{lemma}
\label{le:Bockstein}
$\beta(a)=b$ in $H^{*,*}(DQ_n;\Z/2)$.
\end{lemma}

\begin{proof}
For brevity write $DQ=DQ_n$.
We look at the long exact sequence
\[ \cdots \ra
 H^{1,1}(DQ;\Z) \llra{\times 2} H^{1,1}(DQ;\Z) \ra H^{1,1}(DQ;\Z/2)
\llra{\delta} H^{2,1}(DQ;\Z) \ra \cdots.
\]
The localization sequence (\ref{eq:purity}) for integral cohomology,
together with the identification of $j_!$ in Lemma \ref{lem:Gysin},
shows that $DQ\inc \P^n$ induces an isomorphism on
$H^{1,1}(\blank;\Z)$.  It follows that if $a$ were the mod $2$
reduction of an integral class, it would also be the image of a class
in $H^{1,1}(\P^n;\Z/2)$.  But $* \inc \P^n$ induces an isomorphism on
$H^{1,1}(\blank;\Z/2)$, whereas the class $a$ in $H^{1,1}(DQ;\Z/2)$
restricts to zero on the basepoint.  
We conclude that $a$ can't be the mod $2$ reduction of an integral
class, and therefore $\delta(a)$ is nonzero.

The sequence (\ref{eq:purity}) (again with our knowledge of $j_!$)
also shows 
that $H^{2,1}(DQ;\Z)$ is isomorphic to $\Z/2$, with the generator being in the
image of $H^{2,1}(\P^{n};\Z) \ra H^{2,1}(DQ;\Z)$.  It follows that
$\delta(a)$ is the unique nonzero element of $H^{2,1}(DQ;\Z)$, and the
mod $2$ reduction of $\delta(a)$ is $b$.
\end{proof}

We need one more lemma before stating the final result.

\begin{lemma}
\label{lem:inclusion}
$H^{2k+1,k+1}(DQ_{2k+2};\Z/2) \map H^{2k+1,k+1}(DQ_{2k+1};\Z/2)$
is injective.
\end{lemma}

\begin{proof}
Consider the diagram
\[
\xymatrixcolsep{1.1pc}
\xymatrix{
 H^{2k,k}(Q_{2k+1};\Z/2) \ar[d] &
    H^{2k+1,k+1}(DQ_{2k+2};\Z/2) \ar[d]\ar[l] & 
    H^{2k+1,k+1}(\P^{2k+2};\Z/2) \ar[d]\ar[l] &\ar[l]\\  
 H^{2k,k}(Q_{2k};\Z/2)  &
    H^{2k+1,k+1}(DQ_{2k+1};\Z/2) \ar[l] & 
    H^{2k+1,k+1}(\P^{2k+1};\Z/2) \ar[l] &\ar[l]_-{f}
} 
\]
in which the rows are portions of localization sequences.  
From Proposition \ref{pr:free} and Properties (D) and 
(F), $f$ must be isomorphic to a map
$\M^{1,1}_2 \map \M^{1,1}_2$,  but Lemma \ref{lem:Gysin} tells us 
that this map is the zero map.
Also, the right-most vertical map is an isomorphism.  A
diagram chase would give us the desired result, if we knew that the
left vertical map was injective.  This map is equal to the map
$\CH^k(Q_{2k+1}) \map \CH^k(Q_{2k})$ after reducing modulo 2.
We look at the diagram
\[ \xymatrix{
\Z\ar@{=}[r] &\CH^k(Q_{2k+1}) \ar[r] & \CH^k(Q_{2k}) \ar@{=}[r] &
              \Z\oplus  \Z \\
\Z\ar@{=}[r] & \CH^k(\P^{2k+2}) \ar[r]^\cong \ar[u]^\cong &
\CH^k(\P^{2k+1})\ar[u] \ar@{=}[r]& \Z
}
\]
Lemma \ref{lem:Gysin}
shows that the left vertical map is an isomorphism.
Lemma~\ref{le:evenpullback} identifies the right vertical map as the
diagonal, and from that information the result follows at once.
\end{proof}

\begin{thm} 
\label{th:hyperbolic}
Let $F$ be a field with $\chara(F)\neq 2$.
\begin{enumerate}[(a)]
\item If $n=2k+1$ then 
$H^{*,*}(DQ_n;\Z/2)\iso \M_2[a,b]/(a^2=\rho a+\tau
b,b^{k+1})$
where $a$ has degree $(1,1)$ and $b$ has degree $(2,1)$.
\item If $n=2k$, there exists an element $\epsilon$ in $\M_2^{1,1}$ such
that $H^{*,*}(DQ_n;\Z/2)\iso \M_2[a,b]/(a^2=\rho a+\tau
b,b^{k+1},ab^k=\epsilon b^k)$
where $a$ and $b$ are as in (a).
\end{enumerate}
\end{thm}

\begin{remark}
We haven't been able to identify the class $\epsilon$ in any
nontrivial case.  This is not important for proving the Hopf
condition, but it would be satisfying to resolve the issue
of whether $\epsilon$ is equal to $0$, or $\rho$, or some other element.
\end{remark}

\begin{proof}
For convenience we will drop subscripts and superscripts: $Q=Q_{n-1}$,
$\P=\P^{n}$, and $DQ=DQ_n$.  We know $H^{*,*}(DQ;\Z/2)$ additively, we
just need to determine the ring structure.

Note that the map $H^{2i,i}(\P) \map H^{2i,i}(DQ)$ is surjective
because it is the map $\CH^i(\P) \map \CH^i(DQ)$.  Therefore,
the nonzero element $t$ of $H^{2,1}(\P;\Z/2)$ goes to
the nonzero element $b$ of $H^{2,1}(DQ;\Z/2)$.  Then $t^i$ maps to
$b^i$, and surjectivity implies that $b^i$ must be the unique nonzero
element in $H^{2i,i}(DQ;\Z/2)$ for $1 \leq i \leq \frac{n}{2}$.

Lemma~\ref{le:Bockstein} showed that $\beta(a)=b$.  Since $\beta^2 =
0$ one has $\beta(b) = 0$, so Property (H) implies that
$\beta(ab^i)=b^{i+1}$.  In particular $ab^i$ is nonzero for $0\leq i
\leq \frac{n}{2}-1$.

Now $\tH^{1,1}(DQ;\Z/2)\iso \M_2^{0,0}a$ and $H^{2i-1,i}(DQ;\Z/2)\iso
\M_2^{0,0} \oplus \M_2^{1,1}b^{i-1}$ for $1\leq i \leq \frac{n+1}{2}$,
where the first factor arises from the generator in degree $(2i-1,i)$.
Property (H) and the fact that $\M_2^{2,1}=0$ implies that
$\beta(x)=0$ for any $x\in\M_2^{1,1}b^{i-1}$.  So we cannot have $ab^i
\in \M_2^{1,1}b^{i-1}$.  

Based on our knowledge of $H^{*,*}(DQ;\Z/2)$
as an $\M_2$-module, we can now conclude that when $n=2k$ the classes
$1,b,b^2,\ldots,b^k$ and $a,ab,ab^2,\ldots,ab^{k-1}$ are a free basis
for $H^{*,*}(DQ;\Z/2)$ over $\M_2$.

The argument is slightly harder when $n=2k+1$, because we must show
that $ab^k$ is nonzero (even though its Bockstein is zero).  However,
we already know that $ab^k$ is nonzero in $H^{*,*}(DQ_{n+1};\Z/2)$.
The map $H^{2k+1,k+1}(DQ_{n+1};\Z/2) \ra H^{2k+1,k+1}(DQ_n;\Z/2)$ is an
injection by Lemma \ref{lem:inclusion} and takes $ab^k$ to $ab^k$
by Proposition \ref{pr:DQmap}.  
It follows that $1,b,\ldots,b^k,a,ab,\ldots,ab^k$ is a
free basis for $H^{*,*}(DQ;\Z/2)$ when $n=2k+1$.
 
\medskip

We next identify $a^2$.  This part of the argument exactly parallels
\cite[pp.~22-23]{V2}.  The class $a^2\in \tH^{2,2}(DQ;\Z/2)$ must be a
linear combination over $\M_2$ of the elements $a$ and $b$:
$a^2=Aa+Bb$ where $A\in \M_2^{1,1}$ and $B\in \M_2^{0,1}\iso \Z/2$.
To identify $A$ it's sufficient to look at the image of $a^2$ under
$H^{*,*}(DQ;\Z/2) \ra H^{*,*}(DQ_1;\Z/2)$, since $Aa + Bb$ goes to
$Aa$ under this map by Proposition \ref{pr:DQmap} and the fact that $b
= 0$ in $H^{*,*}(DQ_1;\Z/2)$.  Note that $DQ_1$ is isomorphic to
$\A^1-0$, and one knows that $H^{*,*}(\A^1-0;\Z/2)\iso
\M_2[a]/(a^2=\rho a)$ by \cite[Lem.~6.8]{V2}.  So $A=\rho$.

To identify $B$, let $K$ be the field consisting of $F$ with a square
root of $-1$ adjoined (unless $F$ already has a square root of $-1$,
in which case $K = F$).  Let $DQ_K$ be the base change of $DQ$ along
the map $\spec K \map \spec F$.  Under the induced map
$H^{*,*}(DQ;\Z/2) \map H^{*,*}(DQ_K;\Z/2)$, $\rho$ maps to zero; so
$\rho a + B b$ maps to $B b$.  Hence it suffices to assume that $F$
contains a square root of $-1$ and show that $a^2 =
\tau b$.

Under the map $H^{*,*}(DQ;\Z/2) \ra H^{*}_{et}(DQ;\mu_2^{\tens *})$
the element $\tau$ becomes invertible (cf. Property (I)), 
and so we
can write $a=\tau a'$ (in $H^{*}_{et}(DQ;\mu_2^{\tens *})$), for some
$a'\in H^{1}_{et}(\pt;\mu_2^0)$.  This group is sheaf cohomology with
coefficients in the constant sheaf $\Z/2$; if $\beta_{et}$ is the
Bockstein on \'etale cohomology induced by $0 \ra \Z/2 \ra \Z/4 \ra
\Z/2\ra 0$ one has that $\beta_{et}(a')=(a')^2$ by a standard property
of the Bockstein on sheaf cohomology (the proof is the same as the
one in topology).  Our remarks in Section~\ref{se:etale} show that the
Bocksteins in motivic and \'etale cohomology are compatible, because
$F$ has a square root of $-1$.  So we now compute that
\begin{myequation}
\label{eq:computation}
a^2=\tau^2 (a')^2=\tau^2 \beta_{et}(a') = \tau\cdot \beta(\tau a')=\tau
\cdot \beta(a)=\tau\cdot b
\end{myequation}
in $H^{*}_{et}(DQ;\mu_2^{\tens *})$.
Note that the third equality uses the analog of Property (H) for
\'etale cohomology, together with the fact that $\beta(\tau)=\rho=0$
(by our assumption on $F$).  

As a consequence of (\ref{eq:computation}), we have in particular that
$a^2$ is nonzero in $H^{*,*}(DQ;\Z/2)[\tau^{-1}]$.  But $a^2 = B b$,
so $B$ must be nonzero.  From the sequence (\ref{eq:square}) we recall
that $\M_2^{0,1}=\{0, \tau \}$, so $B=\tau$.  We
have therefore shown that $a^2=\rho a+\tau b \in H^{2,2}(DQ;\Z/2)$.

This finishes part (a) of the theorem.  For part (b) we just observe
that $ab^k\in H^{2k+1,k+1}(DQ;\Z/2)$, and $H^{2k+1,k+1}(DQ;\Z/2)\iso
\M_2^{1,1}b^k$.  So for some $\epsilon \in\M_2^{1,1}$ we have
$\epsilon b^k=ab^k$.  This finishes part (b).
\end{proof}

\begin{remark}
When $n$ is odd, the cohomology of $\DQ_n$ is the same as the
cohomology of the scheme $(\A^n-0)/\pm 1$, which was essentially
computed by Voevodsky in \cite[Th.~6.10]{V2}.  It seems likely that
these two schemes are $\A^1$-homotopy equivalent, but we haven't
proven this.
\end{remark}

\appendix
\label{app:Chow}

\section{Chow groups of quadrics}

This appendix contains a calculation of the Chow rings of
the quadrics $Q_n$, as well as various pushforward and pullback maps.
This is classical, but the details are useful and we don't know a
suitable reference.  We assume a basic familiarity with the Chow ring;
see \cite{F} or \cite[App.~A]{H}.

Let $\CH_i(X)$ be the Chow group of dimension $i$ cycles
on $X$.  If $Z\inc X$ is a closed subscheme there is an exact sequence
\[
\CH_i(Z) \ra \CH_i(X) \ra \CH_i(X-Z) \ra 0
\]
where the first map is pushforward and the second map is restriction.

If $X\subseteq \P^n$ is a closed subscheme, we let $\cone X\subseteq
\P^{n+1}$ denote the projective cone on $X$.  Let $\cone \colon
\CH_i(X) \ra \CH_{i+1}(\cone X)$ be the map sending a cycle to the
projective cone on the cycle, and recall that this is an isomorphism
for $i\geq 0$.  
Also note that $\CH_0(\cone X)=\Z$ no matter what $X$
is.  Finally, recall that $\CH_i(\A^n)=0$ if $i\neq n$, whereas
$\CH_n(\A^n)=\Z$.

When $X$ is nonsingular one defines $\CH^i(X)=\CH_{\dim X-i}(X)$.

The following discussion is modeled on \cite[13.3]{Sw}.

\subsection{The odd-dimensional case}
%
Consider the quadric $Q_{2k+1}\inc \P^{2k+2}$ defined by
$a_1b_1+\cdots+a_{k+1}b_{k+1}+c^2=0$.  We let $j$ be the
inclusion.

\begin{lemma}
\label{lem:odd-pushforward}
For all $0 \leq i \leq 2k+1$, the Chow group $\CH_i(Q_{2k+1})$ is
isomorphic to $\Z$.  The pushforward map $j_*:\CH_i(Q_{2k+1}) \map
\CH_i(\P^{2k+2})$ is an isomorphism if $0 \leq i \leq k$, and is
multiplication by 2 (as a map $\Z \ra \Z$) if $k+1 \leq i \leq 2k+1$.
\end{lemma}

\begin{proof}
The first claim follows immediately from Proposition \ref{pr:free}
and Property (A).

The proof of the second statement is by induction.  The base case
$Q_1$ is isomorphic to $\P^1$, and $Q_1$ is imbedded in $\P^2$ as a
degree two hypersurface.  So $j_*$ is an isomorphism for $i = 0$ and
is multiplication by 2 for $i=1$.

If $Z$ is the closed subscheme defined by $a_1=0$, we know
$Q_{2k+1}-Z\iso \A^{2k+1}$ and $Z\iso \cone Q_{2k-1}$.  The resulting
localization sequence gives us a diagram
\[ \xymatrix{
\CH_i(\cone Q_{2k-1}) \ar[r]\ar[d] & \CH_i(Q_{2k+1}) \ar[r]\ar[d]^{j_*} &
  \CH_i(\A^{2k+1}) \\
\CH_i(\cone \P^{2k}) \ar[r] & \CH_i(\P^{2k+2})
}
\]
in which the top row is exact.
Since $\cone \P^{2k}$ is isomorphic to $\P^{2k+1}$, the bottom
horizontal arrow is an isomorphism for all $0 \leq i \leq 2k+1$,
and both groups in the bottom row are isomorphic to $\Z$.

For $0 \leq i \leq 2k$, the first two groups in the top row are also
isomorphic to $\Z$.  For $0 \leq i \leq k$, the left vertical arrow is
known by induction to be an isomorphism.  The only possibility is that
the map $j_*$ is an isomorphism in this range.

Now for $k+1 \leq i \leq 2k$, 
the left vertical arrow is known by induction
to be multiplication by 2.
Since the upper left horizontal arrow is a surjection,
the only possibility is that 
the map $j_*$ is multiplication by 2.

Finally, for the case $i = 2k+1$ note that $Q_{2k+1}$ is a degree 2
hypersurface in $\P^{2k+2}$.  Thus, the fundamental class $[Q_{2k+1}]$
maps to twice the generator of $\CH_{2k+1}(\P^{2k+2})$.
\end{proof}

By analyzing the above proof, one can give explicit generators for
$\CH_i(Q_{2k+1})$.  If $0 \leq i \leq k$, the generator is the class
of the cycle determined by setting all coordinates equal to zero
except for $b_1, \ldots, b_{i+1}$.  Note that this cycle is isomorphic
to $\P^i$.  On the other hand, if $k+1 \leq i \leq 2k+1$, then the
generator is the class of the cycle determined by setting $a_1,
\ldots, a_{2k+1-i}$ equal to zero.  Note that this cycle is the
iterated projective cone on $Q_{2i-2k+1}$, and also the intersection
of $Q_{2k+1}$ with a copy of $\P^{i+1}$.

We next want to compute the ring structure on $\CH^*(Q_{2k+1})$ as
well as the pullback map $j^*:\CH_i(\P^{2k+2}) \ra
\CH_{i-1}(Q_{2k+1})$.  It is easier to do the latter first.

\begin{prop}
\label{pr:odd-j^*}
The map $j^*\colon\CH_i(\P^{2k+2}) \ra \CH_{i-1}(Q_{2k+1})$ is
an isomorphism if $k+2 \leq i \leq 2k+2$ and is
multiplication by 2 if $ 1 \leq i \leq k+1$.
\end{prop}

\begin{proof}
The projection formula $j_*(a\cdot j^*b)=(j_*a) \cdot b$ gives us
\[ j_*(j^* [\P^i])=j_*([Q_{2k+1}]\cdot j^*[\P^i])=j_*([Q_{2k+1}])\cdot [\P^i] 
   = 2[\P^{2k+1}]\cdot [\P^i] = 2[\P^{i-1}].
\]
In other words the composition 
$j_* j^*: \CH_i(\P^{2k+2}) \map \CH_{i-1}(\P^{2k+2})$ is
multiplication by 2.  When $1 \leq i \leq k+1$, the map
$j_*$ is an isomorphism, so $j^*$ must be multiplication by 2.
When $k+2 \leq i \leq 2k+2$, the map $j_*$ is multiplication by 2,
so $j^*$ must be an isomorphism.
\end{proof}

It is now easy to deduce the ring structure on $\CH^*(Q_{2k+1})$,
using the map from $\CH^*(\P^{2k+2})$.  Note that when $k=0$ we are
looking at $Q_1\iso \P^1$, and so $\CH^*(Q_1)$ is isomorphic to
$\Z[x]/x^2$, where $x$ has degree 1.

\begin{thm}
\label{th:CHoddquadric}
If $k \geq 0$,
then $\CH^*(Q_{2k+1})\iso \Z[x,y]/(x^{k+1}-2y,y^2)$,
where $x$ has degree $1$ and $y$ has degree $k+1$.
\end{thm}
%

\begin{proof}
The map $j^*\colon\CH^i(\P^{2k+2}) \map \CH^i(Q_{2k+1})$ (which 
now preserves the grading because we are grading by codimension)
is an isomorphism if $0 \leq i \leq k$ and is multiplication
by 2 if $k+1 \leq i \leq 2k+1$.  This follows immediately from the
previous proposition simply by regrading.

Let $t$ be the generator of $\CH^1(\P^{2k+2})$, and let $x=j^*(t)$.
Then $x^{k+1}=j^*(t^{k+1})$ is twice a generator of
$\CH^{k+1}(Q_{2k+1})$, and we let $y$ be this generator.  
The desired isomorphism of rings follows immediately from our
knowledge of the groups $\CH^*(Q_{2k+1})$ and the
description of $j^*$ in the previous paragraph.
\end{proof}

\subsection{The even-dimensional case}
\label{se:even-Chow}
This case is a little harder.  The quadric $Q_{2k}$ is defined by
$a_1b_1+\cdots+a_{k+1}b_{k+1}=0$.  As before, let $j$ be the inclusion
$Q_{2k} \map \P^{2k+1}$.  Many of the results from the previous
section carry over to this section with identical proofs.

The base case is $Q_0$, which is $\pt \amalg \pt$.
Note that $j_*:\CH_0(Q_0) \map \CH_0(\P^1)$ is the fold map
$\Z \oplus \Z \map \Z$.

We already know from Proposition \ref{pr:free} that the Chow group
$\CH_i(Q_{2k})$ is isomorphic to $\Z$ for all $0 \leq i \leq 2k$,
except that $\CH_n(Q_{2k})$ is isomorphic to $\Z \oplus \Z$.  The same
arguments as in the proof of Lemma \ref{lem:odd-pushforward} allow us
to conclude that the pushforward map $j_*:\CH_i(Q_{2k}) \map
\CH_i(\P^{2k+1})$ is an isomorphism if $0 \leq i \leq k-1$, is
multiplication by 2 if $k+1 \leq i \leq 2k$, and is the fold map if
$i=k$.  We summarize these facts (with cohomological grading) 
in the following lemma, which we
state because it is critical for the computations in Section
\ref{se:computation}.

\begin{lemma}
\label{lem:Gysin}
For any $n$, the map
$j_*\colon \CH^i(Q_{n-1}) \ra \CH^{i+1}(\P^{n})$ is multiplication by 2
for $0 \leq i < \frac{n-1}{2}$, and is an isomorphism for 
$\frac{n-1}{2} < i \leq n-1$.  
If $n$ is odd, then it is the fold map
$\Z\oplus \Z \ra \Z$ for $i=\frac{n-1}{2}$.  
\end{lemma}

Once again, one can give explicit generators for $\CH_i(Q_{2k})$.  For
$i \neq k$, the description of these generators is the same as in the
odd case.  For $i = k$, one generator is determined by $b_1 = b_2 =
\cdots = b_{k+1} = 0$, and the other generator is determined by $a_1 =
b_2 = \cdots = b_{k+1} = 0$.  We let $\alpha$ and $\beta$ represent
these two codimension $k$ cycles.

\begin{lemma}
Let $\alpha'$ be the cycle determined
by $a_1 = a_2 = \cdots = a_{k+1} = 0$, and let
$\beta'$ be the cycle determined by
$b_1 = a_2 = \cdots = a_{k+1} = 0$.
If $k$ is odd, then $\alpha = \alpha'$ and $\beta = \beta'$ in $\CH_n(Q_{2k})$.
If $k$ is even, then $\alpha = \beta'$ and $\beta = \alpha'$ in $\CH_n(Q_{2k})$.
\end{lemma}

\begin{proof}
If $k$ is odd, then $\alpha$ and $\alpha'$ do not meet.  Then
\cite[Th.~XII.4.III]{HP} says that $\alpha$ and $\alpha'$ are
rationally equivalent.
Similarly, $\beta = \beta'$.

The same argument applies to the even case.
\end{proof}

\begin{lemma}
\label{lem:intersection}
Let $[*]$ be the fundamental class of a point in $\CH^{2k}(Q_{2k})$.
If $k$ is odd, then $\alpha \cdot \alpha = 0 = \beta \cdot \beta$
and $\alpha \cdot \beta = [*]$ in the Chow ring $\CH^*(Q_{2k})$.
If $k$ is even, then $\alpha \cdot \alpha = [*] = \beta \cdot \beta$
and $\alpha \cdot \beta = 0$.
\end{lemma}

\begin{proof}
When $k$ is odd, $\alpha \cdot \alpha = \alpha \cdot \alpha'$.
However, $\alpha$ and $\alpha'$ do not intersect, 
so $\alpha \cdot \alpha' = 0.$  Similarly, $\beta \cdot \beta = 0$.
Now $\alpha$ and $\beta'$ intersect transversely at a point, so
$\alpha \cdot \beta = \alpha \cdot \beta' = [*]$.
Similar arguments apply to the even case.
\end{proof}

As in Proposition \ref{pr:odd-j^*}, the map $j^*\colon \CH_i(\P^{2k+1}) \ra
\CH_{i-1}(Q_{2k})$ is an isomorphism if $k+2 \leq i \leq 2k+1$ and is
multiplication by $2$ if $1 \leq i \leq k$.  After regrading by
codimension, this says $j^*\colon
\CH^i(\P^{2k+1}) \ra \CH^i(Q_{2k})$ is an isomorphism for $0\leq i <k$
and multiplication by $2$ for $k<i\leq 2k$.  The same argument with
the projection formula also shows that when $i = k$, $j^*$ takes the
generator to $u\alpha+(2-u)\beta$ for some $u\in \Z$.

\begin{lemma}
\label{le:evenpullback}
The map $j^*\colon \CH^k(\P^{2k+1}) \map \CH^k(Q_{2k})$
sends the generator $t^k$ to $\alpha + \beta$.
\end{lemma}

\begin{proof}
We already know that $j^*(t^{2k}) = 2[*]$, where
$[*]$ is the fundamental class of a point in $Q_{2k}$ and is also
the generator of $\CH^{2k}(Q_{2k})$.  Therefore,
\[
2[*] = j^*(t^{2k}) = (j^*(t^k))^2 = (u\alpha + (2-u)\beta)^2 =
u^2 \alpha^2 + 2u(2-u)\alpha\beta + (2-u)^2 \beta^2.
\]
If $k$ is odd, Lemma \ref{lem:intersection} lets us rewrite this
equation as
$2[*] = 2u(2-u)[*]$, so $u = 1$.
If $k$ is even Lemma \ref{lem:intersection} gives
$2[*] = (u^2 + (2-u)^2)[*]$, so again $u = 1$.
\end{proof}

\begin{thm}
\label{th:CHevenquadric}
If $k$ is odd, then there is an isomorphism of rings
$\CH^*(Q_{2k})\iso \Z[x,y]/(x^{k+1} - 2xy, y^2)$, where $x$ has degree
$1$ and $y$ has degree $k$.  
If $k$ is even, then $\CH^*(Q_{2k})\iso \Z[x,y]/(x^{k+1} - 2xy,
x^{k+1}y, y^2-x^ky)$, where $x$ has degree $1$ and $y$ has degree $k$.
\end{thm}

\begin{proof}
Let $t$ be the generator $[\P^{2k}]$ of $\CH^1(\P^{2k+1})$, and let
$x=j^*(t)$.  Then $j^*(t^i) = x^i$.  As we know that $j^*$ takes
generators to generators for $0 \leq i \leq k-1$, it follows that
$x^i$ is a generator for $\CH^i(Q_{2k})$ in these dimensions.

Now $x^k = j^*(t^{k}) = \alpha + \beta$ by the previous lemma.
If we let $y$ equal $\alpha$, then $x^k$ and $y$ are 
two generators for $\CH^k(Q_{2k})$.  Note that
Lemma~\ref{lem:intersection} implies that $x^ky=[*]$
since $\alpha (\alpha + \beta) = [*]$ in both the even and odd cases.

Next we can compute that
\[ j_*(x^i \cdot y) = 
  j_*(j^*(t^i) \cdot y)=t^i\cdot j_*y=t^i\cdot t^k=t^{k+i}\]
for $1 \leq t \leq k$.
Since $j_*$ is an isomorphism in codimension $k+i$,
it follows that $x^iy$ is a generator in
$\CH^{k+i}(Q_{2k})$.  

Now $x^{k+1}=j^*(t^{k+1}) = j^*(j_*(xy))=2xy$.  Also, for dimension
reasons $x^{k+1}y = 0$.  Finally, Lemma \ref{lem:intersection} shows
that $y^2 = 0$ if $k$ is odd and $y^2 = [*]= x^ky$ if $k$ is even.

Thus, we have shown that the additive generators for $\CH^*(Q_{2k})$
are $1$, $x$, $x^2, \ldots, x^k$, $y$, $xy,\ldots, x^ky$, where the elements
are listed in order of increasing degree.  Moreover, we have
constructed a ring map from the desired ring to $\CH^*(Q_{2k+1})$,
which is an additive isomorphism.  
\end{proof}


\bibliographystyle{amsalpha}

\end{document}